\documentclass[smallextended]{svjour3}       
\usepackage{latexsym}
\usepackage{amssymb}
\usepackage{amsmath}
\usepackage[mathscr]{eucal}
\usepackage{graphicx}
\usepackage{hyperref}
\usepackage{caption}
\usepackage{subcaption}

\newtheorem{Lemma}{Lemma}

\newtheorem{Theorem}{Theorem}

\renewcommand{\qed}{\hfill{\ \ \rule{2mm}{2mm}} \vspace{0.2in}}

\renewcommand{\thefigure}{\arabic{figure}}
\begin{document}



\title{Tracking the Mean of a Piecewise Stationary Sequence}
\author{ \textbf{Ghurumuruhan Ganesan}
\thanks{E-Mail: \texttt{gganesan82@gmail.com} } \\
\ \\
IISER, Bhopal}
\date{}
\maketitle

\begin{abstract}
In this paper we study the problem of tracking the mean of a piecewise stationary sequence of independent random variables. First we consider the case where the transition times are known and show that a direct running average performs the tracking in short time and with high accuracy. We then use a single valued weighted running average with a tunable parameter for the case when transition times are unknown and establish deviation bounds for the tracking accuracy. Our result has applications in choosing the optimal rewards for the multiarmed bandit scenario.

\vspace{0.1in} \noindent \textbf{Key words:} Piecewise stationary sequence, mean tracking, weighted running average, multiarmed bandit.

\vspace{0.1in} \noindent \textbf{AMS 2000 Subject Classification:} Primary: 60C05, 68T10;
\end{abstract}

\bigskip

\renewcommand{\theequation}{\thesection.\arabic{equation}}
\setcounter{equation}{0}
\section{Introduction} \label{intro}

The multiarmed bandit problem (Auer et al. 2002) is important from both theoretical and application perspectives and has implications in many fields including economics, machine and financial markets, to name a few. The basic setup is that there are multiple arms each equipped with a reward sequence and the average reward in each arm is different. The user (bandit) wants to latch on to the best arm in as short time as possible and therefore must balance exploration together with exploitation to minimize its overall regret. There are many known bounds regarding the exploration exploitation tradeoff (Lai and Robbins (1985), Burnetas and Katehakis (1996)) and for a recent survey of both theoretical and computational aspects, we refer to Slivkins (2019).

Recently (Besbes et al. 2014), there has been interest in the case when the reward sequence of an arm displays nonstationarity. In this case, a direct running average may not be feasible anymore and it is of interest to establish methods for detecting the best arm in a dynamic manner. Auer et al. (2018) propose dynamic tracking for the two armed bandit problem and proceed by comparing the empirical means of both the arms (obtained from sliding window running averages) and declaring that a transition has occurred when the difference is significant. For more details and theoretical bounds on the regret we refer to Lattimore and Szepesv\'ari (2020).

In general, to compute dynamic running averages and invoke deviation estimates, we would  have to store at least~\(O(\log{t})\) samples for~\(t\) rounds.  In this paper, we describe and analyse a \emph{single} valued recursive estimator that uses only knowledge of the minimum epoch durations and automatically performs tracking when the means change.

The paper is organized as follows: In Section~\ref{sec_learn}, we describe the problem of tracking the mean of a piecewise stationary sequence and show that an unweighted running average performs well if the transition times are known. We also define a single valued weighted running average for the case when the transition times are known and describe its tracking performance in our main result Theorem~\ref{thm_seq}. Next in Section~\ref{sec_proof}, we prove Theorem~\ref{thm_seq}.

\renewcommand{\theequation}{\thesection.\arabic{equation}}
\setcounter{equation}{0}
\section{Piecewise Stationary Sequence} \label{sec_learn}
For integer~\(t \geq 1,\) let~\(\{X_j\}_{1 \leq j \leq t}\) be positive independent random variables with~\(0 \leq X_j \leq 1\) and suppose there are deterministic integers~\(M= M(t)\) and~\(1 = s_1 < s_2 < \ldots < s_{M} = t\) such that
\begin{equation}\label{piecewise}
\mathbb{E}X_j=m(k) \text{ for each } s_{k} \leq j < s_{k+1} \text{ and } 1 \leq k \leq M-1.
\end{equation}
In other words, the random variables~\(\{X_j\}\) are not necessarily stationary but enjoy ``epochs" of stationarity in the sense of~(\ref{piecewise}). Throughout this paper, we consider only first order stationarity, i.e., stationarity in the mean and for concreteness, define~\(s_k \leq j <s_{k+1}\) to be the~\(k^{th}\) epoch and~\(\{s_k\}\) to be \emph{transition times}.

It is of interest to ``track" the mean~\(\{m(k)\}\) as it varies from epoch to epoch, using a \emph{randomly} sampled subsequence of~\(\{X_j\}\) as follows. Let~\(\{Z_j\}_{1 \leq j \leq t}\) be independent Bernoulli random variables (that are also independent of~\(\{X_{j}\}\)) with
\begin{equation}\label{samp_prob}
\mathbb{P}(Z_j=1) = \epsilon_j = 1-\mathbb{P}(Z_j=0).
\end{equation}
For~\(t \geq 1,\) let~\({\cal E}(t) = \{1 \leq j \leq t: Z_j=1\}\) be the set of all sampling times and let~\(T(1) < T(2) < \cdots < T(w)\) be the random indices in~\({\cal E}(t).\) Using~\(\{X_{T(l)}\}\) we are interested in tracking the mean sequence~\(\{m(k)\}.\)

In the context of armed bandit sequential learning,~\(\{X_j\}\) denotes the reward sequence of an  arm and~\(\{Z_j\}\) denotes the exploration sequence used to study and gain information regarding the reward sequence. There are usually multiple arms each with its own reward sequence and after every transition, the goal is \emph{track} the arm with the largest reward, i.e., to latch on to the best arm as quickly as possible, using the rewards observed in all arms during the exploration sequence. 

If the transition times~\(\{s_k\}\) are \emph{known}, then we can track~\(\{m(k)\}\) simply using law of large numbers, provided each epoch in~(\ref{piecewise}) is of large enough duration. We use the following deviation estimates to demonstrate our argument.
\begin{Lemma}\label{lemmax}
\((a)\) Let~\(\{W_j\}_{1 \leq j \leq r}\) be independent Bernoulli random variables with~\[\mathbb{P}(W_j = 1) = 1-\mathbb{P}(W_j = 0) > 0.\] Setting~\(S_r := \sum_{j=1}^{r} W_j,\) we have for~\(0 < \epsilon \leq \frac{1}{2}\) that
\begin{equation}\label{conc_est_f}
\mathbb{P}\left(\left|S_r - \mathbb{E}S_r\right| \geq \epsilon \mathbb{E}S_r\right) \leq 2\exp\left(-\frac{\epsilon^2}{4}\mathbb{E}S_r\right).
\end{equation}
\((b)\) Let~\(\{U_j\}_{1 \leq j \leq r}\) be positive independent random variables satisfying~\(0 \leq U_j \leq 1\) and set~\(V_r := \sum_{j=1}^{r}\lambda_jU_j\) where~\(\lambda_j >0\) are positive numbers. For any~\(\epsilon > 0\) we have that
\begin{equation}\label{conc_mom_f}
\mathbb{P}\left(\left|V_r - \mathbb{E}V_r \right| \geq \epsilon \mathbb{E}V_r \right) \leq 2\exp\left(-\frac{\epsilon^2(\mathbb{E}V_r)^2}{\sum_{j=1}^{r}\lambda_j^2}\right).
\end{equation}
\end{Lemma}
For a proof of~(\ref{conc_est_f}) and~(\ref{conc_mom_f}), we refer to Appendix~\(A\) of Alon and Spencer (2008).

Indeed, suppose that the transition times~\(\{s_k\}\) are known and satisfy 
\begin{equation}\label{piece_cond_ax}
\min_{k}(s_{k+1}-s_k) \geq t^{\gamma_0} \text{ and } \min_{k} m(k) \geq \mu_0,
\end{equation}
for some constants~\(0 < \gamma_0 \leq 1\) and~\(\mu_0 > 0.\) Set~\(\epsilon_j = \frac{1}{j^{\beta}}\) where~\(0 < \beta < 1\) is a constant. For~\( 0 < \gamma < \gamma_0,\) let~\({\cal E}_{k,\gamma}(t) \subset {\cal E}(t)\) be the set of all sampling indices lying between~\(s_k\) and~\(s_k+t^{\gamma}\) and
define
\begin{equation}\label{rl_def}
R_l := \frac{1}{\#{\cal E}_{k,\gamma}(t)} \sum_{T(j) \in {\cal E}_{k,\gamma}(t)}X_{T(j)}
\end{equation}
for each~\(l\) satisfying~\(s_k+t^{\gamma} \leq T(l) < s_{k+1}.\) For~\(s_k \leq T(l) < s_k+t^{\gamma},\) we set~\(R_l  = 0.\) In words, we collect all observations lying between~\(s_k\) and~\(s_k + t^{\gamma}\) and compute a single average as defined in~(\ref{rl_def}), whose right hand side does not depend on~\(l.\)

Using~(\ref{conc_est_f}) and~(\ref{conc_mom_f}), we obtain the deviation estimate for~\(R_l\) as follows. Since~\(s_k+t^{\gamma} < t\) for each~\(k,\) the expected size of~\({\cal E}_{k,\gamma}(t)\) is
\begin{equation}
\mathbb{E}\left(\#{\cal E}_{k,\gamma}(t)\right) = \sum_{j=s_k}^{s_k+t^{\gamma}} \frac{1}{j^{\beta}} \geq \sum_{j=t}^{t+t^{\gamma}} \frac{1}{j^{\beta}}  \geq \frac{t^{\gamma}}{t^{\beta}}. \label{ek_size}
\end{equation}
where~\(\#A\) denotes the cardinality of~\(A.\)  Therefore using~(\ref{ek_size}) together with~(\ref{conc_est_f})  we  get that
\begin{equation}\label{ek_dev2}
\mathbb{P}\left(\#{\cal E}_{k,\gamma}(t) \geq C_1t^{\gamma-\beta}\right) \geq 1-\exp\left(-C_2t^{\gamma-\beta}\right)
\end{equation}
for some constants~\(C_1,C_2>0,\) not depending on~\(k.\) Similarly~\(\mathbb{E}(\#{\cal E}(t)) = \sum_{j=1}^{t} \frac{1}{j^{\beta}} \leq \frac{2t^{1-\beta}}{1-\beta}\) and so again using~(\ref{conc_est_f}), we get that
\begin{equation}\label{ek_dev_two}
\mathbb{P}\left(\#{\cal E}(t) \leq \frac{4t^{1-\beta}}{1-\beta}\right) \geq 1-\exp\left(-C_3t^{1-\beta}\right)
\end{equation}
for some constant~\(C_3 > 0.\) There are in total~\(M \leq t\) epochs and so defining~\[E_{nice} := \bigcap_{k=1}^{M-1}\left\{\#{\cal E}_{k,\gamma}(t) \geq C_1t^{\gamma-\beta}\right\} \bigcap \left\{\#{\cal E}(t) \leq \frac{4t^{1-\beta}}{1-\beta}\right\}\] we get from~(\ref{ek_dev2}),~(\ref{ek_dev_two}) and the union bound that
\begin{equation}\label{f_nice_est}
\mathbb{P}(E_{nice}) \geq 1- t\exp\left(-C_2t^{\gamma-\beta}\right) - \exp\left(-C_3t^{1-\beta}\right).
\end{equation}

Given a realization~\(\omega \in E_{nice},\) the expected value~\(\mathbb{E}_{\omega}(R_l)  = m(k)\)  and so the Azuma-Hoeffding inequality~(\ref{conc_mom_f}) together with the second condition in~(\ref{piece_cond_ax}) imply that
\begin{align}\label{rl_est}
\mathbb{P}_{\omega}\left(\left|R_l-m(k)\right| \geq \epsilon m(k)\right) &\leq 2\exp\left(-C_4\epsilon^2 \#{\cal E}_{k,\gamma}(t) \right) \nonumber\\
&\leq 2\exp\left(-C_4C_1\epsilon^2t^{\gamma-\beta}\right)
\end{align}
for each~\(s_k+t^{\gamma} \leq T(l) < s_{k+1}\) and some constant~\(C_4 = C_4(\mu_0)>0\) not depending on~\(k.\) There are~\(M \leq t\) epochs and so defining~\(F_{good}\) to be the intersection of the event that
\begin{equation}\label{et_est}
\#{\cal E}(t) \leq \frac{4t^{1-\beta}}{1-\beta}
\end{equation}
and the event that
\begin{equation}\label{track_def}
|R_l-m(k)| \leq \frac{m(k)}{t^{b}} \text{ for each }k \text{ and each } s_k+t^{\gamma} \leq T(l) < s_{k+1},
\end{equation}
we use~(\ref{ek_dev_two}) and~(\ref{rl_est}) with~\(\epsilon = \frac{1}{t^{b}}\) together with the union bound, to get that~\(F_{good}\) occurs with probability at least
\begin{equation}\label{f_good_est}
1-t\exp\left(-C_2t^{\gamma-\beta}\right) - \exp\left(-C_3t^{1-\beta}\right)-2t\exp\left(-C_4C_1 t^{\gamma(1-\beta)-2b} \right).
\end{equation}
The expression in~(\ref{f_good_est}) converges to one as~\(t \rightarrow \infty\) provided we fix~\(b < \frac{\gamma(1-\beta)}{2}.\)



The relations~(\ref{et_est}) and~(\ref{track_def}) capture the essential requirements we seek in tracking:\\
\((1)\) From~(\ref{et_est}), we get that the number of samples used for tracking is much smaller than the overall duration~\(t.\)\\
\((2)\) The condition~(\ref{track_def}) states that once our tracking scheme comes close to the mean in an epoch, it stays that way for the rest of the epoch, even when new samples are encountered.\\
\((3)\) The term~\(t^{\gamma}\) could be interpreted as an upper bound for the time needed to track the mean with an accuracy of~\(1-\frac{1}{t^b}.\) Since we have chosen~\(\gamma < \gamma_0\) strictly, the tracking time in any epoch is much smaller than the total epoch duration, by~(\ref{piece_cond_ax}).


The caveat in the above analysis is that we have crucially used the knowledge of the transition times~\(\{s_k\}\) for the tracking. This in fact allowed us to construct~\(R_l\) by collecting only the samples observed between~\(s_k\) and~\(s_k + t^{\gamma}.\) In practice,~\(\{s_k\}\) are unknown and so we are interested in ``universal" schemes that do not depend on~\(\{s_k\}.\) In our main result below, we use the following recursive estimator to perform ``automatic" tracking.   For~\(l \geq 1\) let~\(\alpha_l := 1-\frac{1}{l^{\delta}}\) where~\(0 < \delta < 1\) is a constant so that~\(\alpha_1=0.\) Setting~\(Y_0 := 0,\) we define the~\(l^{th}\) weighted running average for~\(l \geq 1\) as
\begin{equation}\label{rec_eq}
Y_{l} := \alpha_l Y_{l-1} + (1-\alpha_l) X_{T(l)}.
\end{equation}
The term~\(\delta\) could be interpreted as the importance or weight given to the past observations in deriving the current average. 

We have the following result regarding the tracking performance of~\(Y_l.\)
\begin{Theorem}\label{thm_seq}
Suppose~(\ref{piece_cond_ax}) holds and we set~\(\epsilon_j = \frac{1}{j^{\beta}}.\) For every~\(0 < \beta < \gamma < \gamma_0\) and every~\(0 < \delta < \frac{\gamma-\beta}{1-\beta},\) there are constants~\(b,\theta > 0\) such that
\begin{equation}\label{e_good_est_ax}
\mathbb{P}\left(E_{good}\right) \geq 1- \exp\left(-t^{\theta}\right)
\end{equation}
where~\(E_{good} = E_{good}(b,\delta,\gamma)\) is event that(\ref{et_est}) and~(\ref{track_def}) hold, with~\(R_l\) replaced by~\(Y_l.\)
\end{Theorem}
In words, knowing only~\(\gamma_0\) and fixing the tracking tolerance~\(\gamma,\) we can choose our parameters~\(\beta\) and~\(\delta\) such that~\(Y_l\) tracks the mean~\(\{m(k)\}\) with high probability.


\renewcommand{\theequation}{\thesection.\arabic{equation}}
\setcounter{equation}{0}
\section{Proof of Theorem~\ref{thm_seq}} \label{sec_proof}
Let~\(\omega\) be a realization of the sequence~\((Z_1,\ldots,Z_t)\) so that~\(\{T(j)\}_{1\leq j \leq L}\) are fixed. We let~\(\mathbb{P}_{\omega}\) be the conditional distribution of~\(\mathbb{P},\) given the sequence~\(\omega.\) Also, for~\(1 \leq k \leq M-1\) let~\[s_k + t^{\gamma}   \leq T(L_k) < T(L_k+1) < \cdots < T(U_k) \leq s_{k+1}-1\] be the set of random indices in~\({\cal E}(t)\) lying between~\(s_k + t^{\gamma}\) and~\(s_{k+1}-1.\)

We prove Theorem~\ref{thm_seq} in three steps. In the first step, we obtain estimates for~\(L_k,U_k\) and their difference. Using these estimates, we show in the second step that the expected value of the running average is close to the current expected value. Finally, we use the Azuma-Hoeffding inequality~(\ref{conc_mom_f}) in Lemma~\ref{lemmax}, to show that the true running averages are close to their respective expected values. Details follow.

\emph{Step 1}: If~\({\cal E}(r_1,r_2) =\{ r_1 \leq j \leq r_2 : Z_j=1\}\) then~\(\#{\cal E}(r_1,r_2) = \sum_{j=r_1}^{r_2} Z_j\) and so
\[\mathbb{E}\left(\#{\cal E}(r_1,r_2)\right) =  \sum_{j=r_1}^{r_2} \frac{1}{j^{\beta}} := a(r_1,r_2)\] where
\begin{equation}\label{a_bounds}
\frac{r_2-r_1}{r_2^{\beta}} \leq a(r_1,r_2) \leq 1+ \int_{r_1}^{r_2} \frac{dx}{x^{\beta}} \leq 1+\frac{r_2^{1-\beta}-r_1^{1-\beta}}{1-\beta}
\end{equation}
for all~\(r\) large. In terms of the notation following~(\ref{samp_prob}), we have~\({\cal E}(r) = {\cal E}(1,r)\) and from~(\ref{a_bounds}), we see that~\(a(1,r) \geq \frac{r-1}{r^{\beta}} \geq \frac{r^{1-\beta}}{2}\) and also that~\(a(1,r) \leq 1+\frac{r^{1-\beta}}{1-\beta} \leq \frac{2r^{1-\beta}}{1-\beta}\) for all~\(r\) large. Consequently, setting~\[E_{samp}(r) := \left\{a(1,r)\left(1-\frac{1}{\log{r}}\right) \leq \#{\cal E}(r) \leq a(1,r)\left(1+\frac{1}{\log{r}}\right)\right\} \] and using the deviation estimate~(\ref{conc_est_f}) with~\(\epsilon = \frac{1}{\log{r}},\) we have that
\[\mathbb{P}\left(E_{samp}(r)\right) \geq 1-\exp\left(-2C_1\frac{r^{1-\beta}}{(\log{r})^2}\right),\] for some constant~\(C_1 > 0.\) We now define
\[E_{exp} := \bigcap_{r \geq (\log{t})^{2/(1-\beta)}} E_{samp}(r)\] and get from the union bound that
\begin{equation}
\mathbb{P}(E_{exp}) \geq 1-\sum_{r \geq (\log{t})^{2/(1-\beta)}} \exp\left(-2C_1\frac{r^{1-\beta}}{(\log{r})^2}\right) \label{e_exp_est_two}
\end{equation}
for all~\(t\) large.

For~\(0 < \epsilon <1-\beta\) strictly and all~\(r \geq (\log{t})^{2/(1-\beta)} =: a\) with~\(t\) large, we have
\[\exp\left(-2C_1\frac{r^{1-\beta}}{(\log{r})^2}\right) \leq \exp\left(-r^{\epsilon}\right)\] and so from~(\ref{e_exp_est_two}), we get
\[\mathbb{P}(E_{exp}) \geq 1- \sum_{r \geq a}\exp\left(-r^{\epsilon}\right) = 1-\sum_{l\geq 0}\exp\left(-(l+a)^{\epsilon}\right).\]
Using the mean value estimate~\((l+a)^{\epsilon} - a^{\epsilon} \geq \epsilon l a^{\epsilon-1} =: \lambda l\) we further get that
\begin{align}
\mathbb{P}(E_{exp}) &\geq 1-\exp\left(-a^{\epsilon}\right) \sum_{l \geq 0}\exp\left(-l\lambda\right)  \nonumber\\
&\geq 1-\exp\left(-a^{\epsilon}\right)\int_0^{\infty} \exp\left(-\lambda (x-1)\right)dx \nonumber\\
&= 1-\frac{1}{\lambda} \exp\left(-a^{\epsilon}\left(1-\frac{\epsilon}{a}\right)\right) \nonumber\\
&\geq 1- \epsilon a^{1-\epsilon}\exp\left(-\frac{a^{\epsilon}}{2}\right) \nonumber
\end{align}
since~\(\epsilon < 1\) and~\(a > 2\) for all~\(t\) large. Setting~\(\epsilon\) to be slightly larger than~\(\frac{3(1-\beta)}{4}\) we then get that
\begin{equation}\label{e_exp_est}
\mathbb{P}(E_{exp}) \geq 1- \exp\left(-(\log{t})^{3/2}\right)
\end{equation}
for all~\(t\) large. We assume henceforth that~\(E_{exp}\) occurs and this obtains upper and lower bounds for each~\(L_k\) and~\(U_k.\) 

For a lower bound on the difference~\(L_{k+1}-U_k,\) we argue as in~(\ref{ek_size}) to get that
\begin{equation}\label{lk_bound}
\mathbb{E}(L_{k+1}-U_k) \geq \sum_{j=s_{k+1}}^{s_{k+1}+t^{\gamma}} \frac{1}{j^{\beta}}  \geq  t^{\gamma - \beta}.
\end{equation}
Choosing~\(\gamma > \beta\) strictly and setting~\[E_{diff} := \bigcap_{k=1}^{M-1} \left\{ L_{k+1}-U_k \geq \frac{t^{\gamma-\beta}}{2}\right\},\] we get from the deviation estimate~(\ref{conc_est_f}) and union bound that
\begin{equation}\label{e_diff_est}
\mathbb{P}(E_{diff}) \geq 1- t\exp\left(-C_2t^{\gamma-\beta}\right)
\end{equation}
for some constant~\(C_2 > 0,\) since the number of transitions~\(M \leq t.\)  Setting~\(E_{tot} := E_{exp} \cap E_{diff}\) we get from~(\ref{e_exp_est}),~(\ref{e_diff_est}) and the union bound that
\begin{equation}\label{e_tot_est}
\mathbb{P}(E_{tot}) \geq 1- \exp\left(-(\log{t})^{3/2}\right)-t\exp\left(-C_2t^{\gamma-\beta}\right).
\end{equation}

We assume henceforth that~\(E_{tot}\) occurs so that for all~\(r \geq (\log{t})^{2/(1-\beta)}\) and~\(t\) large, we have
\begin{equation}\label{er_bounds}
\frac{r^{1-\beta}}{4} \leq \#{\cal E}(r) \leq  \frac{4r^{1-\beta}}{1-\beta},
\end{equation}
by~(\ref{a_bounds}).




\emph{Step 2}: Let~\(m(k)\) be the value of the mean~\(\mu_{j}\) for~\(s_k \leq j < s_{k+1}.\) In this step, we show that for each realization~\(\omega \in E_{tot},\) the expected value of the running average~\(\nu_{l}:=\mathbb{E}_{\omega}Y_{l}\) is close to the current mean~\(m(k)\) for all~\(l\) satisfying~\( L_k \leq l \leq U_k.\)

Defining~\(W_j := X_{T(j)},\)  we see that each~\(Y_{l}\) is a linear sum of~\(\{W_j\}\) since by expanding the recursion~(\ref{rec_eq}) we get
\begin{equation}\label{eq_one}
Y_{l} = \sum_{j=1}^{l}W_j \theta_{j,l}
\end{equation}
where~\[\theta_{j,l}= (1-\alpha_{j})\prod_{i=j+1}^{l}\alpha_i\] for~\(1 \leq j \leq l-1\) and~\(\theta_{l,l}=1-\alpha_l.\) By definition~\(\alpha_1=0\) and~\(\sum_{j=1}^{l}\theta_{j,l}=1\) and so for each~\(L_1 \leq l \leq U_1\) we have that~\(\nu_l = m(1).\)

Next we consider the case~\(L_2 \leq l < U_2\) for which we have
\begin{eqnarray}
  \nu_l &=& m(1)\sum_{j=L_0}^{U_1}\theta_{j,l} + m(2)\sum_{j=U_1+1}^{l}\theta_{j,l} \nonumber\\
   &=&  m(1)\prod_{j=U_1+1}^{l}\theta_{j,l} + m(2)\left(1-\prod_{j=U_1+1}^{l}\theta_{j,l}\right). \label{nu_two}
\end{eqnarray}
Since~\(l \geq L_2,\) we have that
\begin{eqnarray}
\prod_{j=U_1+1}^{l}\alpha_j  &\leq& \prod_{j=U_1+1}^{L_2} \alpha_j \nonumber\\
   &=& \prod_{j=U_1+1}^{L_2} \left(1-\frac{1}{j^{\delta}}\right) \nonumber\\
   &\leq& \exp\left(-\sum_{j=U_1+1}^{L_2} \frac{1}{j^{\delta}}\right)  \nonumber\\
   &\leq& \exp\left(-\frac{(L_2-U_1)}{L_2^{\delta}}\right).\label{temp_pf}
\end{eqnarray}

Since~\(\omega \in E_{exp},\) we get that~\(L_k \leq \frac{4t^{1-\beta}}{1-\beta}\) for any~\(k\) and since~\(\omega \in E_{diff}\) as well, we also know that~\(L_{k+1}-U_k \geq \frac{t^{\gamma-\beta}}{2}.\) Using these estimates in~(\ref{temp_pf}), we get that
\begin{equation}\label{alph_exp}
\prod_{j=U_1+1}^{l}\alpha_j \leq \exp\left(-C_3 t^{\gamma-\beta-\delta(1-\beta)} \right)
\end{equation}
for some constant~\(C_3 > 0.\) Choosing~\(0 < \delta < \frac{\gamma-\beta}{1-\beta}\) strictly, we then get from~(\ref{nu_two}) that
\[\left| \nu_{l} - m(2)\right| \leq (m(1)+m(2))\exp\left(-C_3 t^{\gamma-\beta-\delta(1-\beta)} \right)\]
for all~\(L_2 \leq l \leq U_2.\) Continuing this way iteratively, we get that for each~\(L_k \leq l \leq U_k\) that
\begin{eqnarray}\label{mean_conc}
\left| \nu_{l} - m(k)\right| &\leq& \exp\left(-C_3 t^{\gamma-\beta-\delta(1-\beta)} \right)\sum_{j=1}^{M}m(j) \nonumber\\
&\leq& C_4t\exp\left(-C_3 t^{\gamma-\beta-\delta(1-\beta)} \right)
\end{eqnarray}
for some constant~\(C_4 > 0,\) since~\(0 \leq X_j \leq 1\) by Theorem statement and so~\(m(j)\leq 1.\)


\emph{Step 3}: In this step, we use the deviation estimate~(\ref{conc_mom_f}) to show that~\(Y_l\) is close to~\(\nu_l.\)  We therefore begin with an upper bound for~\(\sum_{j=1}^{l}\theta_{j,l}^2.\)

By definition~\(\theta_{l,l}^2 = (1-\alpha_l)^2 = \frac{1}{l^{2\delta}}\) and for~\(1 \leq j \leq l-1\) we have that
\begin{eqnarray}
\theta_{j,l}^2 &=& (1-\alpha_j)^2\prod_{j+1 \leq i \leq l} \alpha_i^2 \nonumber\\
&=& \frac{1}{j^{2\delta}} \prod_{j+1 \leq i \leq l} \left(1-\frac{1}{i^{\delta}}\right)^2 \nonumber\\
&\leq& \frac{1}{j^{2\delta}} \exp\left(-2\sum_{j+1 \leq i \leq l} \frac{1}{i^{\delta}}\right) \nonumber\\
&\leq& \frac{1}{j^{2\delta}} \exp\left(- \frac{2(l-j)}{l^{\delta}}\right). \nonumber
\end{eqnarray}
If~\(j \leq l - l^{\delta}(\log{l})^2,\) then~\( \theta_{j,l}^2 \leq e^{-2(\log{l})^2}\) and if~\(j \geq l - l^{\delta} (\log{l})^2 > \frac{l}{2},\)
then~\(\theta_{j,l}^2 \leq \frac{4}{l^{2\delta}}.\) Combining we then get that
\begin{equation}\label{theta_est}
\sum_{j=1}^{l} \theta_{j,l}^2 \leq l e^{-2(\log{l})^2} + l^{\delta} (\log{l})^2 \cdot \frac{4}{l^{2\delta}} \leq \frac{5(\log{l})^2}{l^{\delta}}
\end{equation}
for all~\(l \geq L_0,\) a large constant.

Plugging~(\ref{theta_est}) into~(\ref{conc_mom_f}) in Lemma~\ref{lemmax} and setting~\(\epsilon = \frac{1}{t^{b}}\) and using the fact that~\(m(k) \geq \mu_0\) a constant, we then get that
\begin{equation}\label{y_conc}
\mathbb{P}_{\omega}\left(|Y_{l} - m(k)| \geq \frac{m(k)}{t^{b}}\right) \leq 2\exp\left(-\frac{5\mu_0^2 l^{\delta}}{t^{2b}(\log{l})^2}\right).
\end{equation}
The estimate~(\ref{y_conc}) is for a single epoch and there are~\(M \leq t\) epochs. Moreover, recalling that~\(\omega \in E_{tot},\) we get from~(\ref{er_bounds}) that~\[ \frac{t^{\gamma(1-\beta)}}{4} \leq L_1 \leq l \leq \frac{4t^{1-\beta}}{1-\beta} \leq t.\] Therefore from~(\ref{y_conc}) and the union bound we get that
\begin{equation}
\mathbb{P}_{\omega}\left(E_{good}\right) \geq 1-2t\exp\left(-\frac{5\mu_0^2 t^{\delta\gamma(1-\beta)-2b}}{(\log{t})^2}\right).
\end{equation}

We choose~\(b\) small so that~\(b < \frac{\delta\gamma(1-\beta)}{2}.\) Averaging over~\(\omega \in E_{tot}\) and using~(\ref{e_tot_est}), we then get~(\ref{e_good_est_ax}). This completes the proof of the Theorem.~\(\qed\)

\subsection*{\em Acknowledgement}
I thank Professors Rahul Roy, Thomas Mountford, C. R. Subramanian and the referee for crucial comments that led to an improvement of the paper. I also thank IMSc and IISER Bhopal for my fellowships.

\subsection*{\em Data Availability Statement}
Data sharing not applicable to this article as no datasets were generated or analysed during the current study.

\subsection*{\em Conflict of Interest}
The authors have no conflicts of interest to declare that are relevant to the content of this article. No funding was received to assist with the preparation of this manuscript.

\bibliographystyle{plain}

\end{document}